\documentclass{llncs}
\usepackage{graphics}

\begin{document}

\title{A numerical method for solution of ordinary differential equations of
       fractional order}

\author{Jacek Leszczynski, Mariusz Ciesielski}

\institute{Technical University of Czestochowa, Institute of Mathematics
           \& Computer Science, ul.~Dabrowskiego~73, 42-200 Czestochowa,
           Poland \\
           \email{jale@k2.pcz.czest.pl, cmariusz@k2.pcz.czest.pl}}

\maketitle

\begin{abstract}
In this paper we propose an algorithm for the numerical solution of arbitrary
differential equations of fractional order. The algorithm is obtained by using
the following decomposition of the differential equation into a~system of
differential equation of integer order connected with inverse forms of
Abel-integral equations. The algorithm is used for solution of the linear and
non-linear equations.
\end{abstract}


\section{Introduction}
In opposite to differential equations of integer order, in which derivatives
depend only on the local behaviour of the function, fractional differential
equations accumulate the whole information of the function in a~weighted form.
This is so called memory effect and has many applications in physics~\cite{Riewe},
chemistry~\cite{Goto}, engineering~\cite{Carpinteri}, etc. For that reason we
need a~method for solving such equations which will be effective, easy-to-use
and applied for the equations in general form. However, known methods used for
solution of the equations have more disadvantages. Analytical methods, described
in detail in~\cite{Oldham,Samko}, do not work in the case of arbitrary
real order. Another analytical method~\cite{Podlubny}, which uses the multivariate
Mittag-Leffler function and generalizes the previous results, can be used only
for linear type of equations. On the other hand, for specific differential
equations with oscillating and periodic solution there are some specific
numerical methods~\cite{r1,r2,r3}. Other numerical methods~\cite{Blank,Diethelm}
allow solution of the equations of arbitrary real order but they work properly
only for relatively simple form of fractional equations.

Let us consider an initial value problem for the fractional differential 
equation
\begin{equation} \label{eq01}
 a_{1}\, _{\tau }D_{t}^{\alpha _{1}}y(t)+
 a_{2}\, _{\tau }D_{t}^{\alpha _{2}}y(t)+
 \ldots +a_{n}\, _{\tau }D_{t}^{\alpha _{n}}y(t)+
 a_{n+1}y(t)=f(t)
\end{equation}
connected with initial conditions
\begin{equation}\label{eq02}
 y^{(k)}(\tau )=b_{k},
\end{equation}
where \( 0<t\leq T<\infty  \), \( a_{i}\in \rm I\! R \), \( a_{1}\neq 0 \),
\( \alpha _{i}\in \rm I\! R_{+} \), \( i=1,\, \ldots ,n \),
\( m_{i}-1\leq \alpha _{i}<m_{i} \), \( m_{i}\in \rm I\! N \),
\( \alpha _{l}>\alpha _{l+1} \) for \( l=1\, ,\ldots ,n-1 \),
\( b_{k}\in \rm I\! R \), \( k=0\, ,\ldots ,m_{1}-1 \), \(
f(t) \) is a~given function defined on the interval \( [0,T] \),
\( y(t) \) is the unknown function which is the solution of eqn.~(\ref{eq01}).

The fractional derivative operator \( D^{\alpha _{i}} \) is defined in the 
Riemann-Liouville sense~\cite{Samko}. We also have other definitions \ of the 
operator like \ Caputo~\cite{Podlubny}, \ Gr\"{u}nwald-Letnikov~\cite{Podlubny},
and Weyl-Marchaud~\cite{Samko}. Regarding to the Riemann-Liouville operator we
can define it as the left-side operator
\begin{equation} \label{eq03}
 (_{a}D^{\alpha }_{t}y)(t)=\frac{1}{\Gamma (m-\alpha )}\frac{d^{m}}{dt^{m}}
 \int\limits ^{t}_{a}\frac{y(\xi )d\xi }{(t-\xi )^{\alpha -m+1}},
\end{equation}
where \( a<t\leq T<\infty  \), \( m-1<\alpha \leq m \), \( m\in \rm I\! N \).
In our study, we also use integral operators defined in the Riemann-Liouville
sense~\cite{Samko}. When \( y(t)\in L_{1}(a,b) \) and \( \alpha >0 \) then
left-side integral operator is defined as
\begin{equation} \label{eq05}
 (_{a}I^{\alpha }_{t}y)(t)=\frac{1}{\Gamma (\alpha )}\int\limits ^{t}_{a}
 \frac{y(\xi )d\xi }{(t-\xi )^{1-\alpha }},\; t>a.
\end{equation}
More informations concerning the operators properties one can find in
literature~\cite{Oldham,Podlubny,Samko}. Applying definitions~(\ref{eq03})
and~(\ref{eq05}) we have the following property
\begin{equation} \label{eq07}
 (_{a}D^{\alpha }_{t}y)(t)=D^{m}(_{a}I^{m-\alpha }_{t}y)(t),
\end{equation}
where \( D^{m} \) represents typical derivative operator of integer order
\( m \). According to~\cite{Podlubny} in the fractional integrals the
composition rule occurs
\begin{equation} \label{eq09}
 _{a}I^{\alpha }_{t}(_{a}I^{\beta }_{t}y(t))=\, 
 _{a}I^{\alpha +\beta }_{t}y(t)=\, 
 _{a}I^{\beta +\alpha }_{t}y(t)=\, 
 _{a}I^{\beta }_{t}(_{a}I^{\alpha }_{t}y(t))\, .
\end{equation}
In the general case the Riemann-Liouville fractional derivatives do not commute
\begin{equation} \label{eq10}
 (_{a}D^{\alpha }_{t}(_{a}D^{\beta }_{t}y))(t)\neq
 (_{a}D^{\alpha +\beta }_{t})(t)\neq 
 (_{a}D^{\beta }_{t}(_{a}D^{\alpha }_{t}y))(t)\, .
\end{equation}
Extending our considerations, the integer operator commutes with the fractional
operator
\begin{equation} \label{eq11}
 (D^{m}_{t}(_{a}D^{\alpha }_{t}y))(t)=(_{a}D^{m+\alpha }_{t})(t)\, ,
\end{equation}
but the opposite property is impossible. The mixed operators
(derivatives and integrals) commute only in the following way
\begin{equation} \label{eq13}
 (_{a}D^{\alpha }_{t}(_{a}I^{\beta }_{t}y))(t)=
 (_{a}D^{\alpha -\beta }_{t})(t)\, ,
\end{equation}
but they do not commute in the opposite way.

We turn our attention to the composition rule in two following ways. First of
all we found in literature~\cite{Riewe,Riewe2} the fact, that authors
neglect the general property of fractional derivatives given by 
eqn.~(\ref{eq10}). Regarding to solution of fractional differential equations we
will apply above properties in the next sections.

\section{Mathematical background}

In this section we concentrate on describing the method which is decomposition
of fractional differential equation into a system of one ordinary differential
equation and inverse forms of Abel-integral equations. We focus on decomposition
ways of arbitrary fractional differential equation~(\ref{eq01}) with initial
conditions~(\ref{eq02}). In our considerations we classify the fractional
differential equations in dependence on the fractional derivative occurrence
as:
\begin{itemize}
 \item one-term equations, in which \( D^{\alpha } \) occurs only one,
 \item multi-term equations, in which \( D^{\alpha } \) occurs much more times.
\end{itemize}
\begin{corollary}
 An arbitrary one-term equation
 \begin{equation} \label{eq15}
  a_{1}\, _{\tau }D_{t}^{\alpha }y(t)+a_{2}y(t)=f(t)\, ,\;
  0\leq m-1\leq\alpha <m,\, m\in \rm I\! N\,,
 \end{equation}
 in which \( a_{1}\neq 0 \) and initial conditions \( y^{(k)}(\tau )=b_{k} \),
 \( (k=0\, ,\ldots ,m-1) \) should be taken into account, can decomposed into the
 following system
 \begin{equation} \label{eq17}
  \left\{ 
   \begin{array}{l}
    D^{m}_{t}z(t)=-\frac{a_{2}}{a_{1}}y(t)+
    \frac{1}{a_{1}}f(t)\\
    y^{(0)}(t)=\sum ^{m-1}_{k=0}b_{k}\frac{(t-\tau )^{k}}{k!}+\, 
    _{\tau }D^{m-\alpha }_{t}z(t)\\
    y^{(1)}(t)=\sum ^{m-1}_{k=1}b_{k}\frac{(t-\tau )^{k}}{k!}+\, 
    _{\tau }D_{t}^{m-\alpha +1}z(t)\\
    \ldots \\
    y^{(m-1)}(t)=b_{m-1}\frac{(t-\tau )^{m-1}}{(m-1)!}+\, 
    _{\tau }D^{2m-\alpha -1}_{t}z(t)
   \end{array}
  \right. .
 \end{equation}
\end{corollary}
\begin{proof}
 We use the property~(\ref{eq07}) and we introduce a~new variable
 \( z(t) \) which we call a~temporary function. For such introduction we have
 $$
   _{\tau }D^{\alpha}_{t}y(t)= D^m_t z(t) \; \; \; \mbox{and} \; \; \; 
   z(t)=(_{\tau}I^{m-\alpha}_t y)(t).
 $$
 The new variable represents the Abel integral equation which solution, found
 in literature~\cite{Samko}, is the left inverse operator
 $$
  y(t)=(_{\tau}D^{m}_{t} z)(t).
 $$
 In dependence on the integer order \( m \) we introduce initial conditions \(
 b_k \) which are multiplied by a~term \( \frac{(t-\tau)^k}{k!} \). The term
 issues from a~kernel of the left-side Riemman-Liouville operator~(\ref{eq03}).
 \qed
\end{proof}
In the system of eqns.~(\ref{eq17}), the first equation is the differential
equation of \( m \) integer order and the next equations represent inverse forms
of Abel-integral equations. We also found in literature~\cite{Blank,Diethelm}
different or similar approaches for solution of the eqn.~(\ref{eq15}).

In mathematical point of view, more interesting is the second class called
multi-term equations. Let us consider the eqn.~(\ref{eq01}) with initial
conditions~(\ref{eq02}). As in previous class we use the rule~(\ref{eq07}).
Following the previous corollary we apply a~system of new variables
\( z_i(t)=(_{\tau}I_t^{m_i-\alpha_i} y)(t)\) and we decompose
eqn.~(\ref{eq01}) into a system of the following equations
\begin{equation}\label{eq20}
 \left\{ 
  \begin{array}{l}
   a_{1}\, D_{t}^{m_{1}}z_{1}(t)+a_{2}\, D_{t}^{m_{2}}z_{2}(t)+
   \ldots +a_{n}\, D_{t}^{m_{n}}z_{n}(t)+\\
   \qquad \qquad +a_{n+1}y(t)=f(t)\\
   z_{1}(t)=(_{\tau }I^{m_{1}-\alpha _{1}}_{t}y)(t)\\
   z_{2}(t)=(_{\tau }I^{m_{2}-\alpha _{2}}_{t}y)(t)\\
   \ldots \\
   z_{n}(t)=(_{\tau }I^{m_{n}-\alpha _{n}}_{t}y)(t)
  \end{array}
 \right. .
\end{equation}
Such system is not a final form and requires additional transformations. In
multi-term class we distinguish two sub-classes depending on
\( \alpha _{i} \) or follows the decomposition on \( m_{i} \):
\begin{itemize}
 \item independent subclass, in which \( m_{i} \) are different derivatives
       orders of integer type \( m_{i}\neq m_{i+1} \) for \( i=1,\ldots ,n-1 \),
 \item dependent subclass, in which one can find some dependence between the
       derivatives order of integer type, e.g. \( m_{1}=m_{2}=m_{3} \) or
       \( m_{1}=m_{4} \), etc.
\end{itemize}

\begin{proposition}
 Taking into consideration the independent subclass it is possible to formulate
 the unique solution of eqn.~(\ref{eq20}) in the following form
 \begin{equation} \label{eq21}
  \left\{ 
   \begin{array}{l}
    D^{m_{1}}_{t}z_{1}(t)=-\frac{1}{a_{1}}\left(a_{2}\,
    D^{m_{2}}_{t}z_{2}(t)+\ldots +a_{n}\,
    D^{m_{n}}_{t}z_{n}(t)\right)+\\
    \qquad \qquad \qquad +\frac{1}{a_1}\left(-a_{n+1}y(t)+
    f(t)\right) \\
    z_{2}(t)=\, _{0}D^{m_{1}-m_{2}-\alpha _{1}+\alpha _{2}}_{t}z_{1}(t)\\
    z_{3}(t)=\, _{0}D^{m_{1}-m_{3}-\alpha _{1}+\alpha _{3}}_{t}z_{1}(t)\\
    \ldots \\
    z_{n}(t)=\, _{0}D^{m_{1}-m_{n}-\alpha _{1}+\alpha _{n}}_{t}z_{1}(t)\\
    y^{(0)}(t)=\sum ^{m_{1}-1}_{k=0}b_{k}\frac{(t-\tau )^{k}}{k!}+\,
    _{0}D^{m_{1}-\alpha _{1}}_{t}z_{1}(t)\\
    y^{(1)}(t)=\sum ^{m_{1}-1}_{k=1}b_{k}\frac{(t-\tau )^{k}}{k!}+\,
    _{0}D_{t}^{m_{1}-\alpha _{1}+1}z_{1}(t)\\
    \ldots \\
    y^{(m_{1}-1)}(t)=b_{m_{1}-1}\frac{(t-\tau )^{m_{1}-1}}{(m_{1}-1)!}+\,
    _{0}D^{2m_{1}-\alpha _{1}-1}_{t}z_{1}(t)
   \end{array}
  \right. .
 \end{equation}
\end{proposition}

\begin{proof}
 In this proof we try to show a~dependence of temporary functions \( z_i(t) \)
 \( (i=2,...,n)\) on the temporary function \( z_1(t)\). We assume, that initial
 conditions \( b_k=0 \). In such way we can easily construct the  proof. We
 observe that the Abel integral operators in~(\ref{eq20}) have the order
 \( 0<m_i-\alpha_i<1\). We consider the inverse form of such operators and for two
 operators we have
 \begin{equation} \label{eq20a}
  _{\tau}D_t^{\beta} z_2(t) = {_{\tau}}D_t^{\gamma} z_1(t),
 \end{equation}
 where \( 0<\beta<1\) and \( 0<\gamma<1\) respectively. Using the
 property~(\ref{eq07}) we obtain
 $$
  D_t^1(_{\tau}I_t^{1-\beta} z_2(t)) = D_t^1(_{\tau}I_t^{1-\gamma} z_1(t)).
 $$
 Neglecting the operator \( D_t^1 \) we can obtain the following formula
 \begin{equation} \label{eq20b}
  _{\tau}I_t^{1-\beta} z_2(t) = {_{\tau}}I_t^{1-\gamma} z_1(t) + C_{1,2},
 \end{equation}
 where \( C_{1,2} \) is an arbitrary constant. Additionally, we multiply the
 operator \( _{\tau}I_t^{\beta} \) for eqn.~(\ref{eq20b}) and applying the
 property~(\ref{eq09}) we obtain
 \begin{equation} \label{eq20c}
  _{\tau}I_t^1 z_2(t) = {_{\tau}}I_t^{1-\gamma+\beta} z_1(t) +
  {_{\tau}}I_t^{\beta} C_{1,2}.
 \end{equation}
 Differentiating both sides by the operator \( D_t^1 \) and applying the
 property~(\ref{eq13}) (see literature~\cite{Oldham,Podlubny}) we obtain the
 final dependence
 \begin{equation} \label{eq20d}
  z_2(t) = {_{\tau}}D_t^{\gamma-\beta} z_1(t) + C_{1,2}\cdot
  \frac{\beta}{\Gamma(1+\beta)\cdot t^{1+\beta}}.
 \end{equation}
 In general way we state that \( z_i(t) \) \( (i=2,...,n) \) depends on
 \( z_1(t) \) together with additional term represented by the constants
 \( C_{1,i} \).
\qed
\end{proof}

\begin{proposition}
 All constants \( C_{1,i}=0\) \( (i=2,...,n) \) in eqn.~(\ref{eq20d}) are
 equaled to zero for \( \tau=0 \).
\end{proposition}

\begin{proof}
 Taking into consideration the formula~(\ref{eq20d}) we can calculate the constants
 $$
  C_{1,i} = \frac{\Gamma(1+\beta)\cdot t^{1+\beta}}{\beta}\cdot
  \left( z_i(t)-{_{\tau}}D_t^{\gamma-\beta} z_1(t)\right).
 $$
 Additionally, putting the initial condition \( t=\tau=0 \) we obtain all
 \( C_{1,i}=0\). \qed
\end{proof}

\begin{remark}
 According to previous considerations, the initial conditions for temporary
 functions \( z_{i}(t) \) satisfy the following dependence
 \begin{equation} \label{eq22}
  z_{i}(\tau = 0)=0,\: i=1,\ldots ,n\, .
 \end{equation}
\end{remark}
In practical point of view given by formula~(\ref{eq22}), we can solve ordinary
differential equation of integer type with variable coefficients under zeros
initial conditions of temporary functions \( z_{i}(t) \).

Let us consider the next subclass of a~fractional differential equation. We
assume, that some integer orders \( m_{i} \) are the same. This may happen in the
case when some \( \alpha _{i} \) and \( \alpha _{i+1} \), \( (i=1,\ldots ,n-1) \)
belong to the same integer values \( m_{i}=m_{i+1} \). We introduce a~parameter
\( r \) which denotes a~number of temporary functions having the same derivative
operator \( D^{m_{1}} \) of integer type. If \( m_{1}=m_{2}=\ldots =m_{r} \)
for different temporary functions \( z_{i}(t) \), \( (i=1,\ldots ,r) \) then it
can express the functions \( z_{i}(t) \), \( (i=2,\ldots ,r) \) through the
\( z_{1}(t) \) function. It may observe a~relationship in eqn.~(\ref{eq21}) 
for \( z_{2}(t) \) and \( z_{1}(t) \) functions respectively \( z_{2}(t)=\, 
_{\tau }D^{m_{1}-m_{2}-\alpha _{1}+\alpha _{2}}_{t}z_{1}(t) \). We assume an 
equalities of the integer orders \( m_{1}=m_{2} \). The previous relationship
becomes \( z_{2}(t)=\, _{\tau }D^{-\alpha _{1}+\alpha _{2}}_{t}z_{1}(t) \).
Following the property \( (_{\tau }D^{-\alpha }_{t}y)(t)=
(_{\tau }I^{\alpha }_{t}y)(t) \) we have \( z_{i}(t)=\,
_{\tau }I^{\alpha _{1}-\alpha _{i}}_{t}z_{1}(t) \), \( (i=2,\ldots ,r) \).
In the ordinary differential equation presented in the system~(\ref{eq21}) we
change all \( z_{i}(t), \) \( (i=2,\ldots ,r) \) which depends on
\( z_{1}(t) \). Moreover, we introduce a~new temporary function in the following
form
\begin{equation} \label{eq25}
 w(t)=z_{1}(t)+\frac{a_{2}}{a_{1}}\, 
 _{\tau }I^{\alpha _{1}-\alpha _{2}}_{t}z_{1}(t)+\ldots +\frac{a_{r}}{a_{1}}\,
 _{\tau }I^{\alpha _{1}-\alpha _{r}}_{t}z_{1}(t).
\end{equation}
We need to find an inverse form of eqn.~(\ref{eq25}) as
\begin{equation} \label{eq26}
 z_{1}(t)=(1+\frac{a_{2}}{a_{1}}\, 
 _{\tau }I^{\alpha _{1}-\alpha _{2}}_{t}+\ldots +\frac{a_{r}}{a_{1}}\, 
 _{\tau }I^{\alpha _{1}-\alpha _{r}}_{t})^{-1}w(t)\, ,
\end{equation}
where \( (1+\frac{a_{2}}{a_{1}}\, _{\tau }I^{\alpha _{1}-\alpha _{2}}_{t}
+\ldots +\frac{a_{r}}{a_{1}}\, _{\tau }I^{\alpha _{1}-\alpha _{r}}_{t})^{-1} \)
denotes the left inverse operator to the operator \( (1+\frac{a_{2}}{a_{1}}\, 
_{\tau }I^{\alpha _{1}-\alpha _{2}}_{t}+\ldots +\frac{a_{r}}{a_{1}}\, 
_{\tau }I^{\alpha _{1}-\alpha _{r}}_{t}) \). We can apply the known Laplace
transform for eqn.~(\ref{eq25}) that to find the left inverse operator in
expression~(\ref{eq26}). The formula~(\ref{eq26}) is not suitable for practical
purposes and computations.
\begin{proposition}
 Connecting expressions~(\ref{eq21}) and (\ref{eq25})\( \div  \)(\ref{eq26}), the
 solution of the second subclass of the fractional differential equation is
 presented as follows:
 \begin{equation} \label{eq27}
  \left\{ 
   \begin{array}{l}
    D^{m_{1}}_{t}w(t)=-\frac{1}{a_{1}}(a_{r+1}\, 
    D^{m_{r+1}}_{t}z_{r+1}(t)+\ldots +\\
    \qquad +a_{n}\, D^{m_{n}}_{t}z_{n}(t)+a_{n+1}\,
    y(t))+\frac{1}{a_{1}}f(t)\\
    z_{1}(t)=(1+\frac{a_{2}}{a_{1}}\, 
    _{0}I^{\alpha _{1}-\alpha _{2}}_{t}+\ldots +\frac{a_{r}}{a_{1}}\, 
    _{0}I^{\alpha _{1}-\alpha _{r}}_{t})^{-1}w(t)\\
    z_{r+1}(t)=\, _{0}D^{m_{1}-m_{r+1}-\alpha _{1}+
    \alpha _{r+1}}_{t}z_{1}(t)\\
    \ldots \\
    z_{n}(t)=\, _{0}D^{m_{1}-m_{n}-\alpha _{1}+\alpha _{n}}_{t}z_{1}(t)\\
    y^{(0)}(t)=\sum ^{m_{1}-1}_{k=0}b_{k}\frac{(t-\tau )^{k}}{k!}+\, 
    _{0}D^{m_{1}-\alpha _{1}}_{t}z_{1}(t)\\
    y^{(1)}(t)=\sum ^{m_{1}-1}_{k=1}b_{k}\frac{(t-\tau )^{k}}{k!}+\, 
    _{0}D_{t}^{m_{1}-\alpha _{1}+1}z_{1}(t)\\
    \ldots \\
    y^{(m_{1}-1)}(t)=b_{m_{1}-1}\frac{(t-\tau )^{m_{1}-1}}{(m_{1}-1)!}+\, 
    _{0}D^{2m_{1}-\alpha _{1}-1}_{t}z_{1}(t)
   \end{array}
  \right. .
 \end{equation}
\end{proposition}
We do not need to proof the above proposition because it was done in the
previous subclass of the multi-term class.

\section{Numerical treatment and examples of calculations}
In this section we propose an explicit numerical scheme. For numerical solution of
the differential equation of integer order we apply one-step Euler's
method~\cite{Palczewski}. There is no limits to extend above approach to the
Runge-Kutta method. According to results presented by Oldham and
Spanier~\cite{Oldham} we use numerical scheme for integral operator as
\begin{equation} \label{eq29}
 _{0}I_{t}^{\alpha }z_{i}=\frac{h^{\alpha }}{2\Gamma (1+\alpha )}
 \left[ z_{0}(i^{\alpha }-(i-1)^{\alpha })+z_{i}+\sum ^{i-1}_{j=1}z_{i-j}\,
 ((j+1)^{\alpha }-(j-1)^{\alpha })
 \right] \, ,
\end{equation}
which is valid for arbitrary \( \alpha >0 \). We also use an algorithm given by
Oldham and Spanier~\cite{Oldham} for numerical differentiation. The algorithm
depends on \( \alpha  \) range. For \( 0\leq \alpha <1 \) we have
\begin{equation} \label{eq30}
 _{0}D_{t}^{\alpha }z_{i}=\frac{h^{-\alpha }}{\Gamma (2-\alpha )}
 \left[ \frac{(1-\alpha )z_{0}}{i^{\alpha }}+
 \sum ^{i-1}_{j=0}(z_{i-j}-z_{i-(j+1)})\,
 ((j+1)^{1-\alpha }-j^{1-\alpha })\right].
\end{equation}
In literature~\cite{Oldham} one may find the algorithm for highest values.

For limited practical solution of eqn.~(\ref{eq26}) we found in literature
Babenko's method~\cite{Babenko}. In general form of eqn.~(\ref{eq26}) we cannot
consider the left inverse operator. For that, we found only a~simplified formula
limited to the two fractional orders \( \alpha_1, \; \alpha_2 \), which binomial
expansion we can show as
\begin{equation} \label{eq32}
 z_{1i}=\sum ^{\infty }_{j=0}(-1)^{j}(\frac{a_{2}}{a_{1}})^{j}
 \cdot I^{(\alpha _{1}-\alpha _{2})j}w_{i}.
\end{equation}

\begin{figure}[h]
 \resizebox*{0.48\textwidth}{!}{\includegraphics{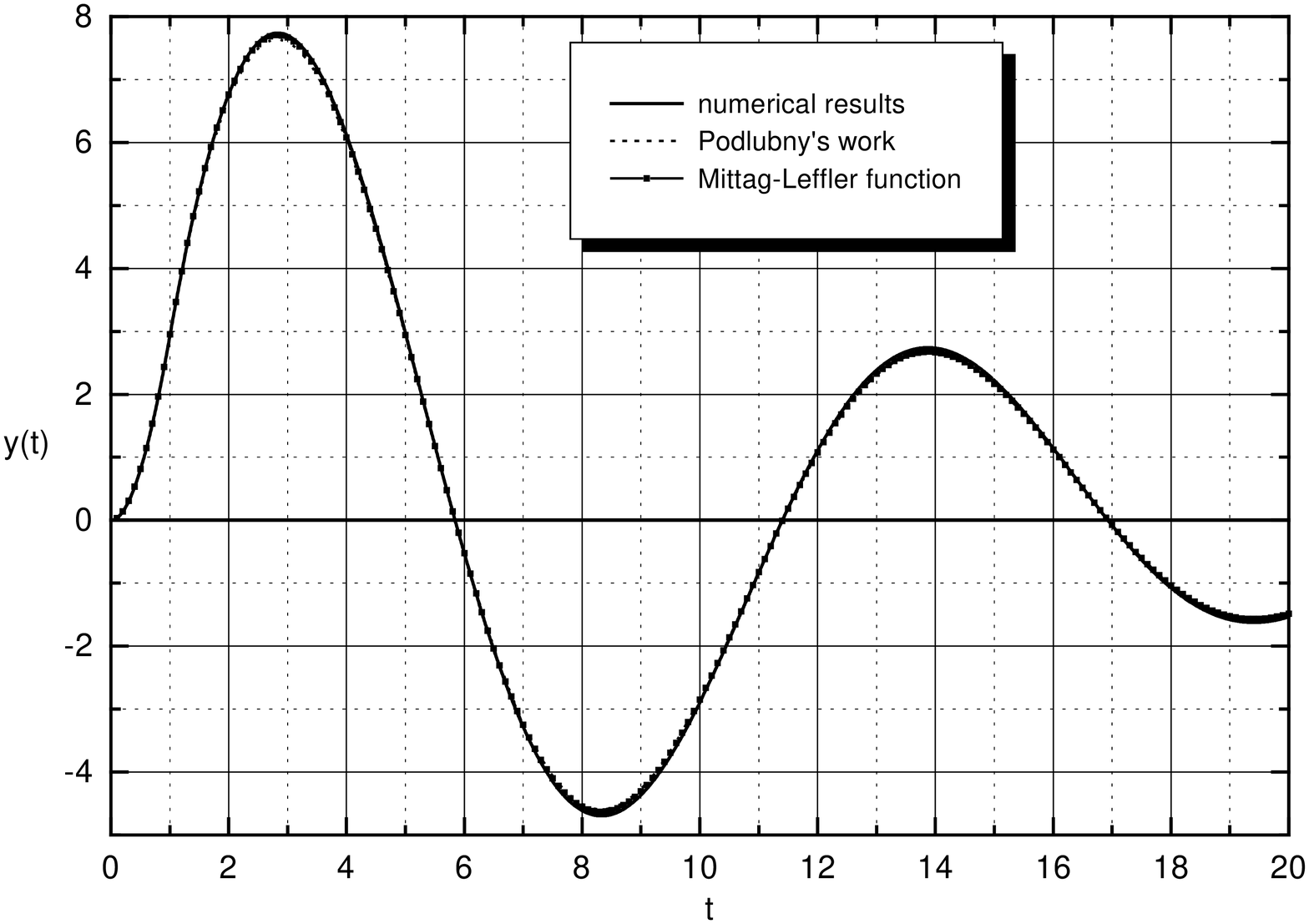}}
 \resizebox*{0.48\textwidth}{!}{\includegraphics{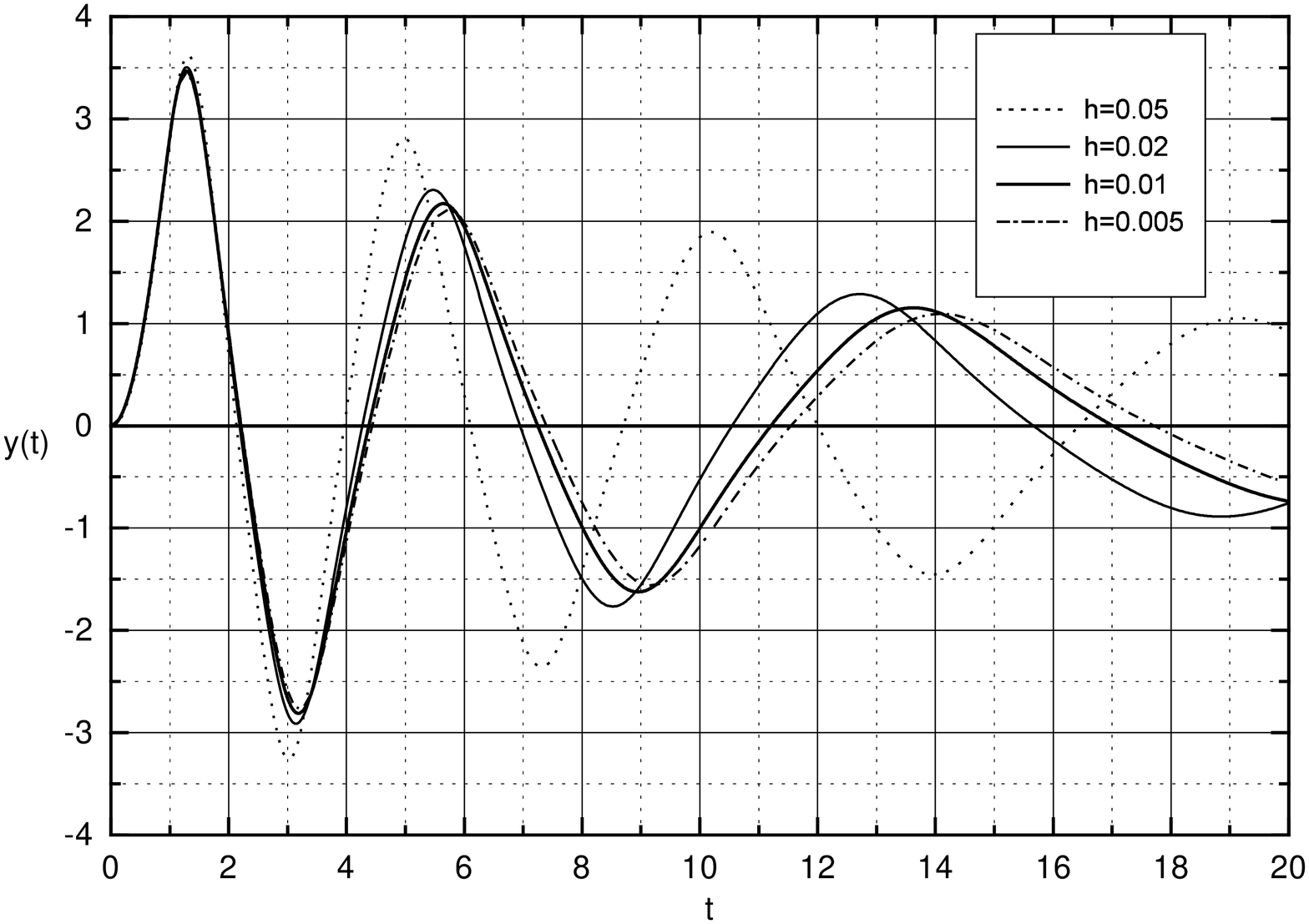}} 
 \caption{\label{fig01} Solution of the Bagley-Torvik equation: \protect\\
          a) comparison of our results with the analytical
             solution~\cite{Podlubny} and with Podlubny's
             work~\cite{Podlubny}, \protect\\
          b) numerical solution of the non-linear form of Bagley-Torvik equation.}
\end{figure}

For some comparison we select the most popular fractional differential equation in
literature~\cite{Podlubny}
\begin{equation} \label{eq44}
 a\cdot _{0}D^{2}_{t}y(t)+b\cdot _{0}D^{1.5}_{t}y(t)+c\cdot y(t)=f(t)
\end{equation}
called the Bagley-Torvik equation. Fig.~\ref{fig01}a presents the qualitative
comparison of the numerical solution generated by our method with the analytical
solution.
We assume: \( a=1 \), \( b=0.5 \), \( c=0.5 \), the function
\( f(t)=\left\{ \begin{array}{l} 8\; \mbox{for} \; 0\leq t<1\\ 0\; \mbox{for} \;
1<t<\infty \end{array}\right.\) and initial conditions \( y(0)=0,\: y'(0)=0 \)
respectively. The analytical solution one may find in \cite{Podlubny}. We also
found a~simple way of solution created in Podlubny's work \cite{Podlubny}.
Fig.~\ref{fig01}a certifies that our numerical results behave good in
comparison to the analytical solution. The great advantage would be the apply of
our approach to the non-linear form of fractional differential equation. Let us
consider the non-linear form of the fractional differential equation
\begin{equation} \label{eq46}
 _{0}D^{2}_{t}y(t)+0.5\cdot \, _{0}D^{1.5}_{t}y(t)+0.5\cdot \, y^{3}(t)=f(t)
\end{equation}
where \( f(t)=\left\{ \begin{array}{l}
8\; \mbox{for} \; 0\leq t<1\\
0\; \mbox{for} \; 1<t<\infty 
\end{array}\right.  \)with initial conditions \( y(0)=0,\: y'(0)=0 \). This is
the Bagley-Torvik equation where non-linear term \( y^3(t) \) is introduced. 
Fig.~\ref{fig01}b shows a~behaviour of the numerical solution of
eqn.~(\ref{eq46}). We can see that the step \( h \) has strong influence to the
solution \( y(t) \). Non-linear fractional differential equations need some
computational tests that to choose right value of the step.

\section{Conclusions}
In this paper we propose a~new method for the numerical solution of arbitrary
differential equations of fractional order. Following to Blank's
results~\cite{Blank} the main advantage of the method is a~decomposition of
fractional differential equation into a~system composed with one ordinary
differential equation of integer order and the left inverse equations of the
Abel-integral operator. We distinguish two classes of such system: one-term
fractional derivative and multi-term fractional derivatives. The comparison
certifies that our method gives quite good results. Summarizing these results,
we can say that the decomposition method in its general form gives a~reasonable
calculations, is the effective method and easy to use and is applied for the
fractional differential equations in general form.

\end{document}